

%

%
%
%
%

\def \reel{{\bf R}}

\def \comp{{\bf C}}

 \baselineskip = 18pt
 \def\n{\noindent}
 \overfullrule = 0pt
 \def\qed{{\hfill{\vrule height7pt width7pt
depth0pt}\par\bigskip}} 
 
 \magnification\magstep1
\magnification\magstep1
\baselineskip = 18pt
\def\n{\noindent}
\def\qed{{\hfill{\vrule height7pt width7pt
depth0pt}\par\bigskip}}

\def\pf{ \medskip \n {\bf Proof.~~}}
\def\ie{{\it i.e.\/}\ }

\centerline{\bf Quadratic forms in unitary operators}\medskip

\centerline{by Gilles Pisier\footnote*{Partially supported by the
N.S.F.}}
\centerline{Texas A\&M University}
\centerline{and} 
\centerline{Universit\'e Paris VI}
\bigskip

Let $u_1,\ldots,u_n$ be unitary operators on a Hilbert space $H$. We will
study the norm
$$\left\|\sum^{i=n}_{i=1} u_i \otimes \bar u_i\right\|\leqno (1)$$
of the operator $\sum u_i \otimes \bar u_i$ acting on the Hilbertian tensor
product $H\otimes_2 \overline H$. Throughout this paper $\overline H$ will
be the complex conjugate of $H$ and $H^*$ the dual space. Of course, we have
canonically $\overline H\simeq H^*$. Therefore, $H\otimes_2 \overline H
\simeq H \otimes _2 H^*$ can be identified with the space $S_2(H)$ of all
Hilbert-Schmidt operators on $H$, equipped with the Hilbert-Schmidt norm,
denoted by $\|~~\|_2$. Then (1) can be rewritten as
$$\left\|\sum^n_{i=1} u_i \otimes \bar u_i\right\| = \sup\left\{ \left\|
\sum^n_1 u_itu^*_i\right\|_2\, \bigg| \, t\in S_2(H),\quad \|t\|_2 \le
1\right\}. \leqno (1)'$$
We will denote by $B(H)$ the space of bounded operators
on $H$ equipped with the usual norm. Note that
$\overline{B(H)}$ can be canonically identified with
$B( \overline H)$. More
generally, let $K$ be another Hilbert space, consider
$a_1,\ldots, a_n \in B(H)$ and $b_1,\ldots, b_n \in
B(K)$, then we can view $\sum a_i\otimes \bar b_i$ as
acting on $H\otimes \overline K \simeq H\otimes K^*$
identified with the Hilbert-Schmidt class $S_2(K,H)$
equipped with the Hilbert-Schmidt norm, again denoted
simply by $\|~~\|_2$.  Then we have $$\left\|\sum^n_1
a_i\otimes \bar b_i\right\| = \sup \left\{\left\| \sum^n_1
a_itb^*_i\right\|_2\, \bigg|\, t\in S_2(K,H), \quad
\|t\|_2 \le 1\right\}.\leqno (2)$$ The left side of (2) is
the norm in the ``minimal'' or ``spatial'' tensor product
of the $C^*$-algebras $B(H)$
 and $\overline{B(K)}$,
which is defined in
full generality as follows
(cf.\ e.g.\
[Ta]): 
 for any Hilbert spaces $H_1,H_2$ and for any $a_k\in
B(H_1)$, $b_k\in B(H_2)$ let us denote by $\sum
a_k\otimes b_k$ the associated linear operator on
$H_1\otimes_2 H_2$ taking $h_1\otimes h_2$ to $\sum
a_k(h_1) \otimes b_k(h_2)$. The norm induced by $B(H_1
\otimes_2 H_2)$ on the algebraic tensor product $B(H_1)
\otimes B(H_2)$ is called the ``minimal'' or ``spatial''
norm. In the sequel, all the norms appearing will be of
this kind, unless specified otherwise.

In matrix notation, of course any element $a\in B(\ell_2)$ can be
represented by a bi-infinite matrix $(a(i,j))$ with complex entries.  The
reader who prefers this framework will recognize that
$a\otimes b$ can be identified with the Kronecker product
of the associated matrices, and $\bar b$ with the matrix
with complex conjugate entries to that of $b$. With this
viewpoint $\sum a_k\otimes \bar b_k$ corresponds to the
matrix $\left(\sum\limits_k a_k(i,j)
\overline{b_k(i',j')}\right)$ where the rows are indexed
by pairs $(i,i')$ and the columns by pairs $(j,j')$.

The expressions appearing in (1) and (2) play a
fundamental r\^ole in the author's recent theory of the
operator Hilbert space $OH$, see [P1].

\n  We
now return to (1). Note that by the triangle inequality we have trivially
$$\left\|\sum^n_1 u_i \otimes \bar u_i\right\|\le n.$$
If $\dim H<\infty$, this cannot be improved and we have
$$\left\|\sum u_i\otimes \bar u_i\right\| =n.\leqno (3)$$
Indeed, $t = Id_H$ is an eigenvector for $t\to \sum u_i
t u^*_i$ associated to the eigenvalue $n$.

More generally, it is easy to see that (3) still holds when $\dim H =
\infty$ if $u_1,\ldots, u_n$ all belong to a finite injective von~Neumann
subalgebra $M\subset B(H)$. However, (3) is not true if we drop the
injectivity assumption, as shown when $M$ is the von~Neumann algebra
(factor actually) associated to the free group $F_n$ on $n$ generators. We
first recall some notation to make this more precise. Let $G$ be any
discrete group (for instance $F_n$). We denote by $\lambda\colon \ G\to
B(\ell_2(G))$ the left regular representation which takes an element
$x$ in $G$ to the unitary operator of left translation by
$x$. Then we denote by $VN(G)$ the von~Neumann algebra
generated by $\lambda(G)$ in $B(\ell_2(G))$.  Now in the
particular case $G=F_n$ let $g_1,\ldots, g_n$ be the
generators of $F_n$ so that
$\lambda(g_1),\ldots,\lambda(g_n)$ are unitary generators
for $VN(F_n)$. Then it is known that $$\left\|\sum^n_1
\lambda(g_i) \otimes \overline{\lambda(g_i)}\right\| =
2\sqrt{n-1} = \left\|\sum^n_1 \lambda(g_i)\right\|. \leqno
(4)$$ Indeed, as we will see below, the left hand side is
the same as $\left\|\sum\limits^n_1 \lambda(g_i)\right\|$
and the latter norm was computed in [AO] and found equal
to the middle side of (4). The results of [AO] were partly
motivated by Kesten's thesis [K], where it is proved that
$$\left\|  \sum^n_1 \lambda(g_i) + \lambda(g_i)^*\right\|
= 2\sqrt{2n-1} \leqno (5)$$ and also that (5) realizes the
minimum of all norms $\left\| \sum\limits_{t\in
S}\lambda(t)\right\|$ when $S$ runs over all possible
symmetric subsets of cardinality $2n$ of any discrete
group $G$.

The next observation, which is our main result, extends Kesten's lower
bound to a more ``abstract'' setting.

\proclaim Theorem 1. For any $n$-tuple $u_1,\ldots, u_n$ of unitary
operators in $B(H)$, we have
$$2\sqrt{n-1} \le \left\|\sum^n_1 u_i \otimes \bar u_i\right\|.\leqno (6)$$
In other words, the right side of (6) is minimal exactly when $u_i =
\lambda(g_i)$.

\pf We will make extensive use of a simple but important
result due to Fell [F], as follows:\ for any (discrete)
group $G$ and unitary representation $\pi\colon\ G\to
B(H)$ the representation $\lambda\otimes \pi$ is unitarily
equivalent to $\lambda\otimes I$ (for a proof see e.g.\
[DCH, p.~469]). We call this ``Fell's absorption
principle''. For convenience, we will apply this to $\bar
\pi$ instead:\ $\lambda\otimes \bar \pi \simeq \lambda
\otimes I$. As a consequence, for any $t_1,\ldots, t_n$ in
$G$ we have $$\left\|\sum^n_1 \lambda(t_i) \otimes
\overline{\pi(t_i)}\right\| = \left\|\sum^n_1
\lambda(t_i)\right\|.\leqno (7)$$ Now, when $G=F_n$ the
data of a unitary representation $\pi\colon \ F_n\to B(H)$
boils down  to the $n$-tuple $u_1,\ldots, u_n$ of the
values of $\pi$ at the free generators $g_1,\ldots, g_n$.
Hence (7) yields that for any choice of unitary operators
$$\left\|\sum^n_1 \lambda(g_i) \otimes \bar u_i\right\| =
\left\|\sum^n_1 \lambda(g_i)\right\|.\leqno (7)'$$
Now by an   inequality  due to Haagerup
(cf. [H1, Lemma 2.4]), we have with the same notation as
in (2) above
$$\left\|\sum\nolimits^n_1 a_i\otimes \bar
b_i\right\| \le \left\|\sum\nolimits^n_1 a_i\otimes \bar
 a_i\right\|^{1/2}
\left\|\sum\nolimits^n_1 b_i
 \otimes \bar b_i\right\|^{1/2}.$$
Therefore, we have 
 $$\left\|\sum
\lambda(g_i)\otimes \bar u_i\right\| \le \left\|\sum
\lambda(g_i) \otimes \overline{\lambda(g_i)}\right\|^{1/2}
\left\|\sum \bar u_i \otimes u_i\right\|^{1/2}.\leqno
(8)$$ Let $s_n = \left\|\sum\limits^n_1
\lambda(g_i)\right\|$. By Fell's principle applied again,
we have $\left\|\sum\limits^n_1 \lambda(g_i) \otimes
\overline{\lambda(g_i)}\right\| =s_n$, hence (7)$'$ and
(8) yield $$s_n \le (s_n)^{1/2} \left\|\sum \bar
u_i\otimes u_i\right\|^{1/2}. $$ Recalling (4)
(and dividing by $(s_n)^{1/2}$), we obtain (6).\qed
\medskip

\n {\bf Remark.} More precisely, the same argument shows that for any
finite subset $S\subset G$ of an arbitrary discrete group $G$ and for any
uniformly bounded representation $\pi\colon \ G\to B(H)$ we have
$$\left\|\sum_{t\in S} \lambda(t)\right\| = \left\|\sum_{t\in S} \lambda(t)
\otimes \overline{\lambda(t)}\right\| \le \sup_{t\in G} \|\pi(t)\|^4
\left\|\sum_{t\in S} \pi(t) \otimes \overline{\pi(t)}\right\|.$$
More generally, for any family $(f(t))_{t\in S}$ of
operators
in $B(H)$, we have
$$\left\|\sum_{t\in S} \lambda(t)
\otimes \overline{\lambda(t)}
\otimes f(t)
\otimes \overline{f(t)}
\right\| \le \sup_{t\in G} \|\pi(t)\|^4
\left\|\sum_{t\in S} \pi(t) \otimes \overline{\pi(t)}
\otimes f(t)
\otimes \overline{f(t)}\right\|.$$

The argument is the same as above but using the version of
the absorption principle given in [DCH, Lemma~2.1, p.~469].

We now apply Theorem~1 to
estimate the constant $c_n$ defined in [JP] for any $n\ge
1$ as follows $$c_n = \inf\left\{\sup_{m\ne m'}
\left\|\sum^{i=n}_{i=1} u^m_i \otimes \overline{u^{m'}_i}
\right\|\right\}$$ where the infimum runs over all
possible choices of {\it infinite\/} sequences
$(u^m_1,\ldots, u^m_n)$ of $n$-tuples of $N_m\times N_m$
unitary matrices with arbitrary size $N_m$. We have
trivially $c_1=1$ and $c_n\le n$ for all $n$. As observed
in [JP], it is true (this is an amusing exercise) that
$c_2=2$, but more importantly (see [JP]) we have $c_n<n$
for any $n\ge 3$. As pointed out by A.~Valette (see [JP,
Remark~2.12] and also Valette's note [Va]), the striking
work of Lubotzky-Phillips-Sarnak (see Lubotzky's book [L])
allows to show that $c_n \le 2\sqrt{n-1}$ for all $n=p+1$
with $p$ prime $\ge 3$. As the next result shows, this
turns out to be an equality, since we have

\proclaim Lemma 2. \quad $2\sqrt{n-1} \le c_n$ for all
$n\ge 1$.

\pf Let $(u^m_i)_{i\le n}$ be a sequence of $n$-tuples with $(u^m_1,\ldots,
u^m_n)$ unitary in the space $M_{N_m}$ of all $N_m\times N_m$ complex
matrices. Let $A$ be the space formed of all families $x = (x_m)_{m\in {\bf
N}}$ with $x_m \in M_{N_m}$ and $\sup\limits_m \|x_m\|_{M_{N_m}} < \infty$.
Equipped with the norm $\|x\| = \sup\|x_m\|_{M_{N_m}}$, $A$ becomes a
$C^*$-algebra. Let ${\cal U}$ be a non-trivial ultrafilter and let $I_{\cal
U}\subset A$ be the (closed two-sided self-adjoint) ideal formed of all
sequences $x = (x_m)_{m\in {\bf N}}$ such that $\lim_{\cal U} \|x_m\|=0$.
Then the quotient space $A/ {I_{\cal U}}$ is a $C^*$-algebra  called the
ultraproduct of $\{M_{N_m}\mid m\in {\bf N}\}$ with respect to ${\cal U}$.
By Gelfand theory we can view $A/ {I_{\cal U}}$ as embedded into
$B(\widehat H)$ for some Hilbert space $\widehat H$.
Let us denote by $\hat u_1,\ldots, \hat u_n$ the unitary elements in
$A/ {I_{\cal U}}$ associated to the families $(u^m_1)_{m\in {\bf
N}},\ldots, (u^m_n)_{m\in {\bf N}}$. We claim that for any $a_1,\ldots,
a_n$ in $B(H)$ (with $H$ arbitrary) we have
$$\left\|\sum \hat u_i\otimes a_i\right\| \le \lim_{m,{\cal U}} \left\|\sum
u^m_i \otimes a_i\right\|.\leqno (9)$$
(Indeed, the quotient mapping $q\colon \ A\to A/ I$ is a
$C^*$-representation, hence $q\otimes I_{B(H)}$ extends to 
a contractive
representation from $A\otimes_{\rm min} B(H)$ to $A/ I \otimes_{\rm min}
B(H)$, see e.g.\ [Pa1] for details.)

\n Now, if we apply (9) with $a_i =  \bar{\hat u}_i \in
B(\overline{\cal H})$, we obtain by Theorem~1
$$\eqalignno{2\sqrt{n-1} \le \left\|\sum \hat u_i \otimes \bar{\hat
u}_i\right\| &\le \lim_{m,{\cal U}} \left\|\sum u^m_i \otimes \bar{\hat
u}_i \right\| = \lim_{m,{\cal U}} \left\|\sum \overline{u^m_i} \otimes \hat
u_i\right\|\cr
&= \lim_{m,{\cal U}} \left\|\sum \hat u_i \otimes
\overline{u^m_i}\right\|\cr
\noalign{\hbox{hence by (9) again}}
&\le \lim_{m,{\cal U}} \lim_{m',{\cal U}} \left\|\sum u^{m'}_i \otimes
\overline{u^m_i}\right\|\cr
\noalign{\hbox{and the last term is of course}}
&\le \sup_{m\ne m'} \left\|\sum^n_{i=1} u^m_i 
\otimes \overline{u^{m'}_i}
\right\|.}$$
Thus we conclude that $2\sqrt{n-1} \le c_n$.\qed

It would be extremely interesting (especially in connection
with Voiculescu's last question in [Vo]) to characterize the
$n$-tuples of unitary operators $(u_1,\ldots, u_n)$ for
which the lower bound in Theorem~1 is attained, i.e.\ for
which $$\left\|\sum^n_1 u_i\otimes \bar u_i\right\| =
2\sqrt{n-1}\qquad (n\ge 3).$$
Although it might be prematurate in view of the lack of examples, we
formulate a conjecture.\medskip

\n  {\bf Conjecture:}\ Let $u_1,\ldots, u_n$ be unitary operators on a
Hilbert space $H$ such that\break $\left\|\sum u_i\otimes \bar u_i\right\|
=
2\sqrt{n-1}$ $(n\ge 3)$. Then the linear mapping which takes $\lambda(g_i)$
to $u_i$ extends to a ``complete contraction'' in the
sense of e.g.\ [Pa1] (actually it might even be
completely isometric). Equivalently this means that there
is a $C^*$-representation $\pi\colon \ VN(F_n)\mapsto B(H)$
and contractive operators $v,w$ in $B(H)$	such that $$u_i
=  v\pi(\lambda(g_i))w \qquad i=1,2,\ldots, n.$$

Note that, by Akemann and Ostrand's characterization of
Leinert sets in [AO], this is true if $u_i=\lambda(x_i)$
with $(x_i)$ any Leinert set
with $n$-elements 
in an arbitrary group $G$. In particular,
if $(u_i)_{i\le n}$ consists of $(\lambda(g_i))_{i\le k}$
and its inverses $(\lambda(g_i)^*)_{i\le k}$ (with $n=2k$)
then the  span of $(u_i)_{i\le 2k}$ is completely
isometric to the span of $(\lambda(g_i))_{i\le 2k}$.

\n However, perhaps this conjecture might only be true or
easier to prove for ``symmetric'' $n$-tuples of the form
$(u_1,u^*_1,u_2,u^*_2,\ldots, u_k,u^*_k)$ with $n=2k$.
Indeed, in this case the conjecture is valid
for group translations:
if $u_i=\lambda(x_i)$
with $(x_i)$ any symmetric set
with $n=2k$-elements 
in an arbitrary group $G$, Kesten [K]
showed that $(x_i)$ must consist of
$k$ free elements and their inverses.

\n {\bf Remark.} More recently, S.\ Szarek (personal communication) found an
alternate proof of (6) closer in spirit to Kesten's proof for group
translations. Let
$$C = \{t\in S_2\mid t\ge 0 \quad \|t\|_2 =1\}.$$
Then (cf.\ e.g.\ [P1, Example~5.6]) for any $u_i$ in $B(H)$
$$\left\|\sum u_i\otimes \bar u_i \right\| = \sup\left\{{tr}\left(\sum
u_itu^*_is\right)\, \Big|\, t,s\in C\right\}.\leqno (10)$$
Note that for any $t,s$ in $C$
$$tr(u_itu^*_is) = tr(s^{1/2}u_i tu^*_is^{1/2})\ge 0.\leqno (11)$$
Moreover, when the family
$(u_1,...,u_n)$ is self-adjoint (\ie when $u_i^*$
also belongs to the family), the
supremum in (10) can be restricted to $t=s$.

\n Let $T =
\sum\limits^n_1 u_i \otimes \bar u_i$ and let $\widetilde
T = \sum\limits^n_{i=1}  \lambda(g_i) \otimes
\overline{\lambda(g_i)}$. Szarek's idea consists in
showing that for any integer $m\ge 1$ and any $t$ in $C$
we have $$\langle (T^*T)^mt,t\rangle \ge \langle
(\widetilde T^*\widetilde T)^m \xi, \xi\rangle\leqno
(12)$$ where $\xi = e\otimes \bar e$ and where $e\in
\ell_2(F_n)$ denotes the basis vector indexed by the unit
element in $F_n$. Note that the normalized trace $\tau$ in
$VN(F_n)$ is given by the formula $$\tau(x) = \langle
x e, e\rangle.\leqno \forall x\in VN(F_n)$$ To verify
(12), note that we can expand $(T^*T)^m$ as a sum of the
form $\sum\limits_{\alpha\in I} u^\alpha \otimes
\overline{u^\alpha}$ where the $u^\alpha$'s are unitaries
of the form $u^*_{i_1}u_{j_1} u^*_{i_2}u_{j_2}\ldots$\ .

\n Now for certain $\alpha$'s, we have $u^\alpha = I$ by
formal cancellation (no matter what the $u_i$'s are), let
us denote by $I'\subset I$ the set of all such $\alpha$'s.
Then by (11) we have for all $t$ in $C$ $$\langle (T^*T)^m
t,t\rangle = \sum_{\alpha\in I} tr(u^\alpha tu^{\alpha
*}t) \ge \sum_{\alpha \in I'} 1 = {\rm card}(I')$$ but by
an elementary counting argument we have $${\rm card}(I') =
\langle(\widetilde T^*\widetilde T)^m \xi,\xi\rangle =
(\tau\otimes\tau)[(\widetilde T^* \widetilde T)^m].$$ Hence
we obtain (12). Therefore $$\|T^*T\| \ge \lim_{m\to \infty}
\langle (T^*T)^m t,t\rangle^{1/m} \ge \lim_{m\to \infty}
((\tau\otimes)\tau[(\widetilde T^* \widetilde T)^m])^{1/m}
= \|\widetilde T^* \widetilde T\|,$$ so that we obtain
$\|T\| \ge \|\widetilde T\|$, whence (6).

The preceding results can be used to give some complementary information
related to the important work of Lubotzky-Phillips-Sarnak 
[LPS] on Ramanujan graphs and distribution of points
on the sphere. To
describe this, we need some notation.

\n Let us denote by $S_{N-1}$ (resp.\ $S^{\comp}_{N-1}$)
the $N$-dimensional sphere in ${\reel}^N$ $({\comp}^N)$
equipped with its standard rotationally invariant (resp.\
unitarily invariant) probability measure. We will denote
simply by $L_2(S_{N-1})$ (resp.\ $L_2(S^{\comp}_{N-1})$)
the associated $L_2$ space and by $L^0_2(S_{N-1})$ (resp.\
$L^0_2(S^{\comp}_{N-1})$) the subspace orthogonal to the
constant function 1.

\n There is a classical unitary representation
$\hat\rho\colon \ SO(N) \mapsto B(L_2(S_{N-1}))$ (called
the ``quasi-regular'' representation) defined by
$$\hat\rho(\omega)f(\cdot) = f(\omega^{-1}(\cdot)),\leqno
\forall \omega \in SO(N) \quad \forall f\in L_2(S_{N-1})$$
and similarly in the complex case. We will denote by
$\rho$ the restriction of $\hat\rho$ to $L^0_2(S_{N-1})$.

Then, Lubotzky-Phillips-Sarnak (see [L]) proved:

\proclaim Theorem 3. {\bf ([LPS]).}
\item{(i)} For any $n$ and any $\omega_1,\ldots, \omega_n$ in $SO(3)$ we
have
$$2\sqrt{n-1} \le \left\|\sum^n_1
\rho(\omega_i)\right\|_{B(L^0_2(S_2))}.$$
\item{(ii)} For any $n$ of the form $n=p+1$ with $p$ prime $\ge 3$, there
are elements $\omega_1,\ldots, \omega_n$ in $SO(3)$ such that
$$\left\|\sum^n_{i=1} \rho(\omega_i)\right\|_{B(L^0_2(S_2))} \le
2\sqrt{n-1}.$$

\n The reader should note that the lower bound (i) is considerably easier
to prove than the upper bound (ii) (the latter uses
Deligne's proof of the Weil conjectures).
An alternate proof of (i) appears
in [CV]. We give another one below. Curiously, both
bounds remain open for $SO(N)$ with $N>3$.

However, we can prove the lower bounds in the complex case, i.e.\ in the
case of $SU(N)$ with $N$ arbitrary. 
We will denote 
by $\rho^{\comp}\colon \ SU(N) \to B(L^0_2(S^{\comp}_{N-1}))$
the quasi-regular representation restricted to the
orthogonal of constant functions. Then we have

\proclaim Theorem 4. Let $N\ge 1$ be arbitrary. Then for any $n$ and any
$\omega_1,\ldots, \omega_n$ in $SU(N)$ we have
$$2\sqrt{n-1} \le
\left\|\sum^{i=n}_{i=1}\rho^{\comp}(\omega_i)\right\|.\leqno(13)$$

Let $\pi\colon \ G\to B(H)$ and $\sigma\colon\ G\to B(H)$
be unitary representations of a group $G$ on a Hilbert
space $H$. Then we denote by $\pi\otimes \bar \sigma\colon
\ G\to B(H \otimes_2 \overline H)$ the unitary
representation defined by $$\pi\otimes \bar \sigma(t) =
\pi(t) \otimes \overline{\sigma(t)}.$$ Consider in
particular the representation $\rho^{\comp}$
on $SU(N)$.
Let us denote by $\{\pi_m\mid m\in {\bf N}\}$ the (finite dimensional) {\it
irreducible\/}  unitary representations which appear in the
decomposition of $\rho^{\comp}$ into irreducible components. By
avoiding repetitions, we may assume that for $m\ne m'$, $\pi_m$ is not
unitarily equivalent to $\pi_{m'}$ (hence the corresponding characters are
orthogonal). Moreover, since $SU(N)$ is compact, all the
representations $\{\pi_m\}$ are {\it finite dimensional\/}.

I am most grateful to Anthony Wassermann for showing me the next result
(probably known to specialists) and the elementary proof below.

\proclaim Lemma 5. Fix $N\ge 1$. Let $(\pi_m)$ be
 associated as above to
$\rho^{\comp} $ on $SU(N)$. Then, one can
extract from it
an infinite subset  $(\sigma_m)$ such that, for each $m\ne
m'$, every irreducible representation appearing in the
decomposition of $\sigma_m \otimes \bar \sigma_{m'}$ is 
included (up to unitary equivalence)
in the original family $\{\pi_m\}$.

\pf  Let $H = L^0_2(S^{\comp}_{N-1})$. Let $m\ge	1$.
Let $H_m \subset H$ be the subspace of all analytic polynomials which are
homogeneous of degree $m$. Then $\rho^{\comp}$ restricted to $H_m$ is
irreducible, let $\sigma_m$ be this representation. We claim that for any $m
\ne m'$, $\sigma_m \otimes \overline{\sigma_{m'}}$ can be written as
$\bigoplus\limits_{\sigma\in\Sigma(m,m')} \sigma$ with $\sigma$ irreducible
subrepresentation of $\rho^{\comp}$.

\n Indeed, consider the linear map $V\colon \ H_m\otimes
\overline{H_{m'}}\to H$ associated to the product, i.e.\
taking $g\otimes\bar h$ to the function $t\to g(t)
\overline{h(t)}$ in $L^0_2(S^{\comp}_{N-1})$. Then  it is not too hard to verify that $V$ is
injective (see below).

\n Moreover, $V$ satisfies
$$V(\sigma_m\otimes\overline{\sigma_{m'}})(\omega) = \rho^{\comp}(\omega)V, \leqno
\forall \omega\in SU(N)$$
in other words $V$ intertwines $\sigma_m\otimes \overline{\sigma_{m'}}$ and
$\rho^{\comp}$ restricted to $V(H_m \otimes \overline{H_{m'}})$. This shows (by
Schur's classical lemma) that every irreducible component $\sigma$ of
$\sigma_m \otimes \overline{\sigma_{m'}}$ appears as a subrepresentation of
$\rho^{\comp}$.

\n We now check the injectivity of $V$
(I am grateful
to E. Straube for showing me this quick argument): 
Let $F(z,w)$ be a polynomial on 
$\comp^N\times \comp^N$,
$m$-homogeneous in $z$, $m'$-homogeneous in $w$
and such that $F(z,\bar z)=0$ for all $z$ in the unit
sphere.
Then, by the homogeneities,  $F(z,\bar z)=0$ for all $z$
in $\comp^N$. This last condition and the analyticity
of $F$ in $\comp^N$ force $F=0$. Indeed,  
  the  derivative $DF$
 (which is $\comp$-linear
by the holomorphy of $F$) must 
satisfy 
$DF(z,\bar z)=0$ for all $z$
in $\comp^N$, and similarly for all successive
derivatives.
Hence the analytic function $F$ must vanish
identically.\qed

\n {\bf Proof of Theorem 4.} By Lemma~2 we have
$$2\sqrt{n-1} \le \sup_{m\ne m'} \left\|\sum^n_{i=1} \sigma_m(\omega_i)
\otimes \overline{\sigma_{m'}(\omega_i})\right\|.$$
Now by Lemma~5 whenever $m\ne m'$ we have $\sigma_m \otimes
\overline{\sigma_{m'}} \simeq \bigoplus\limits_{\sigma\in \sum (m,m')} \sigma$
where $\sum(m,m')$ consists of subrepresentations of $\rho^{\comp}$, hence in
particular we have for all $m\ne m'$
$$\left\|\sum^n_{i=1} \sigma_m (\omega_i) \otimes
\overline{\sigma_{m'}(\omega_i)}\right\| \le
\left\|\sum_1^n \rho^{\comp}(\omega_i)\right\|. $${\qed}

\n {\bf Remark.} The arguments of Valette to show that
$c_n \le 2\sqrt{n-1}$ when $n=p+1$ with $p$ prime $\ge 3$
can be easily described using Theorem~3 (ii) and the
preceding discussion, so we briefly sketch it for the
reader's convenience:\ Let us denote again by $\{\pi_m\}$
the collection of distinct irreducible finite dimensional
unitary representations appearing in $\rho$, but this
time in the real case. The point is that, 
 on $SO(3)$, it
is known that, if $m\ne m'$, all the irreducible components of
$\pi_m\otimes
\overline{\pi_{m'}}$ are subrepresentations of $\rho$ (i.e.\ are in the
family $(\pi_m)$). The  reason behind this
is simply
that {\it all} non-trivial irreducible
representations appear
as subrepresentations of $\rho$ (the latter
fact is no longer true
on 
$SO(N)$ with $N>3$, cf. e.g. [Vi, p. 440-457]). Therefore
we have again $$\left\|\sum^n_{i=1} \pi_m(\omega_i) \otimes
\overline{\pi_{m'}(\omega_i)} \right\| \le \left\|\sum^n_1
\rho(\omega_i)\right\|,$$ so that choosing $u^m_i =
\pi_m(\omega_i)$ we can deduce $$c_n \le \left\|\sum^n_1
\rho(\omega_i)\right\|.\leqno(14)$$
Now by (ii) in Theorem 3, this implies
$c_n \le 2\sqrt{n-1}$. 

\n Note that, by Lemma 2,
(14) also gives a new proof of (i) in Theorem~3.
 
\n However, on
$SO(N)$ with $N>3$, 
 the same argument apparently
does not extend and (i) in Theorem~3  remains open
for $SO(N)$ if $N>3$.

\n {\bf Remark.} Of course the same lower bound  (13)
is valid with the same proof 
for any unitary representation 
$\rho$ on a  group $G$  provided
there
exists an {\it infinite} set $(\sigma_m)$ of {\it finite
dimensional} unitary representations of $G$ such that
every irreducible component
of $\sigma_m\otimes \overline{\sigma_{m'}}$ with
$m\not=m'$ is incuded in $\rho$.

\vskip12pt
\n{\bf Acknowledgement:} This note was conceived
during an extended stay
 at Trinity College and 
 Cambridge University's DPMMS.
 I would like to express my warmest thanks to
    B\'ela Bollob\'as for his 
extremely kind hospitality there.
 \vfill\eject

\centerline{\bf References}

\item{ [AO]} C. Akemann  and P. Ostrand.  Computing
norms in group $C^*$-algebras.  Amer. J. Math. 98
(1976) 1015-1047.

\item{ [CV]} Y. Colin de Verdi\`ere. 
Distribution de points
sur une sph\`ere (d'apr\`es Lubotzky, Phillips et Sarnak).
S\'eminaire N. Bourbaki. Ast\'erisque Soc. Math. France
177-178 (1989) 83-93. 

\item{[DCH]} J. de Canni\`ere   and  U. Haagerup. 
Multipliers of the Fourier algebras of some simple Lie
groups and their discrete subgroups.  Amer. J. Math.   107
(1985), 455-500.

\item{[F]} M. Fell. Weak containment and induced representations of
groups. Canad. J. Math. 14 (1962) 237-268.

 \item{[H1]} U. Haagerup. Injectivity and decomposition of completely
 bounded maps in ``Operator algebras and their connection with Topology and
 Ergodic Theory''. Springer Lecture Notes in Math. 1132 (1985) 91-116.
 
\item{[JP]} M. Junge and G. Pisier. Bilinear forms on
exact operator
spaces and $B(H)\otimes B(H)$. Geometric and functional
analysis (GAFA) 5 (1995) 329-363.

\item{[K]} H. Kesten.  Symmetric random walks on groups.
Trans.  Amer. Math.  Soc. {92} (1959) 336-354.

\item{[L]} A. Lubotzky.
Discrete groups, expanding graphs and invariant measures.
Progress in Math. 125.
Birkhauser, 1994.

\item{[LPS]} A. Lubotzky, R. Phillips and P. Sarnak. Hecke
operators
and distributing
points on $S^2$,I. Comm. Pure and Applied Math. 39 (1986)
S149-186.

\item{[Pa1]} V. Paulsen.  Completely
bounded maps and dilations. Pitman Research Notes 146.
Pitman Longman (Wiley) 1986.

\item{[P1]} G. Pisier. The operator Hilbert space $OH$,
complex interpolation and tensor norms.
Memoirs Amer.
Math. Soc. (to appear)
 
\item{[Ta]} M. Takesaki. Theory of Operator Algebras I.
Springer-Verlag New-York 1979.

\item{[Va]} A. Valette. An application of Ramanujan graphs
to $C^*$-algebra tensor products. To appear.

\item{[Vi]} N. J. Vilenkin. Special functions and the
theory
of group
representations. Translations  Math. Monographs 22.
Amer. Math. Soc. Providence, RI, 1968.

\item{[Vo]} D. Voiculescu. About quasi-diagonal operators.
Integr. Equat. Oper. Th. 17 (1993) 137-149.

\vskip12pt

Texas A\&M University

College Station, TX 77843, U. S. A.

and

Universit\'e Paris VI

Equipe d'Analyse, Case 186,
 
75252 Paris Cedex 05, France

 \end

More generally, for any Hilbert spaces $H_1,H_2$ and for any $a_i\in B(H_1)$,
$b_i\in B(H_2)$ we will denote by $\sum a_k\otimes b_k$ the associated
linear operator on $H_1\otimes_2 H_2$ taking $h_1\otimes h_2$ to $\sum
a_k(h_1) \otimes b_k(h_2)$. The norm induced by $B(H_1 \otimes_2 H_2)$ on
the
algebraic tensor product $B(H_1) \otimes B(H_2)$ is called the ``minimal''
or ``spatial'' norm. In the sequel, all the norms appearing will be of this
kind, unless specified otherwise.

In matrix notation, of course any element $a\in B(\ell_2)$ can be
represented by a bi-infinite matrix $(a(i,j))$ with complex entries.  The
reader who prefers this framework will recognize that
$a\otimes b$ can be identified with the Kronecher product
of the associated matrices, and $\bar b$ with the matrix
with complex conjugate entries to that of $b$. With this
viewpoint $\sum a_k\otimes \bar b_k$ corresponds to the
matrix $\left(\sum\limits_k a_k(i,j)
\overline{b_k(i',j')}\right)$ where the rows are indexed
by pairs $(i,i')$ and the columns by pairs $(j,j')$.

\font\ninebf=msym9

\def\C{{\ninebf C}}
$\rho^{\C }$

$\rho_{\comp}$

$\comp$

\end